\documentclass[10pt]{amsart}
\usepackage[a4paper,margin=27mm]{geometry}
\usepackage[T1]{fontenc}
\usepackage{lmodern,amsmath,amssymb,mathtools,booktabs,microtype,xcolor}
\usepackage[colorlinks=true,linkcolor=blue!50!black,citecolor=blue!50!black,urlcolor=blue!50!black]{hyperref}
\hypersetup{pdftitle={Gopakumar--Vafa invariants and Mordell--Weil lattices},pdfauthor={Lutian Zhao}}
\newtheorem{theorem}{Theorem}[section]
\newtheorem{proposition}[theorem]{Proposition}
\newtheorem{corollary}[theorem]{Corollary}
\newtheorem{conjecture}[theorem]{Conjecture}
\newtheorem{lemma}[theorem]{Lemma}
\theoremstyle{definition}

\newtheorem{example}[theorem]{Example}
\theoremstyle{remark}
\newtheorem{remark}[theorem]{Remark}
\newcommand{\OO}{\mathcal O}
\newcommand{\Q}{\mathbb Q}
\newcommand{\Z}{\mathbb Z}
\newcommand{\C}{\mathbb C}
\newcommand{\PP}{\mathbb P}
\DeclareMathOperator{\Pic}{Pic}
\DeclareMathOperator{\NS}{NS}
\DeclareMathOperator{\MW}{MW}
\DeclareMathOperator{\Ext}{Ext}
\DeclareMathOperator{\Hom}{Hom}
\DeclareMathOperator{\Sym}{Sym}
\DeclareMathOperator{\Chow}{Chow}
\DeclareMathOperator{\IC}{IC}
\DeclareMathOperator{\ch}{ch}
\DeclareMathOperator{\length}{length}
\DeclareMathOperator{\tr}{tr}
\newcommand{\pH}{{}^{p}\!\mathcal H}
\newcommand{\HH}{\mathbb H}
\numberwithin{equation}{section}
\title[GV invariants and Mordell--Weil lattices]{Gopakumar--Vafa invariants and Mordell--Weil lattices}
\author{Lutian Zhao}
\address{Kavli Institute for the Physics and Mathematics of the Universe (WPI), The University of Tokyo, Japan}
\email{lutian.zhao@ipmu.jp}
\date{November 18, 2026}
\subjclass[2020]{14N35, 14J27, 14D20}
\keywords{Gopakumar--Vafa invariants, Mordell--Weil lattice, height pairing, elliptic Calabi--Yau threefold, vanishing cycles}

\begin{document}
\begin{abstract}
Let $X\to B$ be a Weierstrass fibration with integral fibers from a smooth Calabi--Yau threefold to a del Pezzo surface. For a smooth curve $D\subset B$ with smooth inverse image $Y$, we decompose the moduli fiber of stable sheaves with Euler characteristic one and projected cycle $D$. In fiber degree $k$, its components are indexed by all sections $C$ of $Y\to D$ with $\tfrac12\langle C,C\rangle\leq k$, where $\langle\, ,\,\rangle$ is the Mordell--Weil height pairing. For infinitesimally rigid sections we compute the GV contributions, recovering the $E_8$ theta product over a $(-1)$-curve. We formulate the corresponding conjecture for families of sections.
\end{abstract}
\maketitle

\section{Introduction}

Huang, Katz, Klemm, and Wang (HKKW) compute refined BPS numbers for compact elliptic Calabi--Yau threefolds using Wilson loops \cite[\S\S3.3--3.4 and~5.1]{HKKW}. In their calculation, distinct sections can contribute to the same curve class \cite[\S\S5.1.4 and~5.1.7]{HKKW}. Fix a smooth curve $D$ in the base. Which sections over $D$ occur in a given fiber degree, and what does each contribute to the Gopakumar--Vafa (GV) invariants after fibers are added? Since the projected cycle is $D$ with multiplicity one, each support has exactly one horizontal component, of degree one over $D$. Theorem~\ref{thm:fixed} describes the moduli fiber as the section varies, and Theorem~\ref{thm:main} computes the ambient GV contribution of each infinitesimally rigid section.

In Gromov--Witten (GW) theory, such contributions enter the series of virtual counts of stable maps through the predicted GV multicover formula
\begin{equation}\label{eq:GW-GV-overview} \sum_{\substack{\beta>0\\g\geq0}}N^{\mathrm{GW}}_{g,\beta}\lambda^{2g-2}Q^\beta=\sum_{\substack{\beta>0,\ h\geq0\\m\geq1}}\frac{n_{h,\beta}}{m}\bigl(2\sin(m\lambda/2)\bigr)^{2h-2}Q^{m\beta}, \end{equation}
where $\beta>0$ runs over nonzero effective curve classes, $N^{\mathrm{GW}}_{g,\beta}$ is the connected genus-$g$ GW invariant, $n_{h,\beta}$ is the genus-$h$ GV number in class $\beta$, $\lambda$ is the genus variable, and the formal symbols satisfy $Q^\beta Q^\gamma=Q^{\beta+\gamma}$ \cite[Part~II, equation~(3.2)]{GV}. The agreement of these numbers with those defined by vanishing cycles below is conjectural in general \cite{MaulikToda}.

For an elliptic threefold, one can group the terms of this series by base and fiber degrees. Alim--Scheidegger and Klemm--Manschot--Wotschke studied the modular structure and holomorphic anomaly equations that arise in this way \cite{AS,KMW}. At fixed base degree, Huang, Katz, and Klemm predicted that the all-genus partition function is expressed in meromorphic Jacobi forms with prescribed poles. In their formula the fiber parameter is the modular variable and the string coupling is the elliptic variable \cite[equations~(1.1)--(1.2)]{HKK2015}. Their later work extends the calculation to refined invariants and examines the dependence on complex structure and orientation \cite[\S\S3.5--3.7]{HKK2021}.

The Jacobi transformation has a monodromy interpretation. Schimannek computed the action of relative Fourier--Mukai monodromy on brane charges \cite[\S3]{SchMon}, and Pioline and Schimannek relate this action to the transformation of the fixed-base partition function \cite[\S4]{PiSch}. On the enumerative side, Oberdieck and Pixton prove the quasi-Jacobi and holomorphic anomaly equations for GW invariants of rational elliptic surfaces \cite{OP}. In Pandharipande--Thomas (PT) theory one takes virtual counts of stable pairs $\OO_X\to\mathcal F$ with $\mathcal F$ pure of dimension one and zero-dimensional cokernel \cite{PT}. Here Oberdieck and Shen prove the elliptic transformation law for reduced base classes \cite[Theorem~1]{OS}.

HKKW, on the other hand, fix the fiber degree and allow the base degree to vary. Over $\PP^2$ they propose a finite expansion of the compact partition function in local Wilson loops \cite[\S3.4.1, equation~(3.58)]{HKKW}, whose values are computed in \cite{KimKimKim,HLW,WangLoops}. For this expansion, a section and its added fibers count toward the same degree. Our first task is therefore to determine the degree of the section itself. Appendix~\ref{app:family} formulates the conjecture for families motivated by HKKW's expansion.

We work over $\C$. Let $B$ be a smooth projective del Pezzo surface, so that $-K_B$ is ample, and let $\pi:X\to B$ be a Weierstrass fibration with smooth projective total space and zero section $\sigma:B\to X$. Locally over $B$, the fibration is defined in Weierstrass form by
\[ y^2=x^3+ax+b,\qquad a\in H^0(B,K_B^{-4}),\quad b\in H^0(B,K_B^{-6}). \]
We assume that the discriminant divisor $4a^3+27b^2=0$ is reduced with only ordinary cusps at the transverse intersections of $a=0$ and $b=0$, and that all fibers of $\pi$ are integral. We further assume that $X$ has trivial canonical bundle, $H^1(X,\OO_X)=0$, $H_2(X,\Z)$ is torsion-free, and the Picard and curve groups decompose integrally as
\begin{equation}\label{eq:lattices} \Pic(X)=\Z[S]\oplus\pi^*\Pic(B),\qquad H_2(X,\Z)=\sigma_*H_2(B,\Z)\oplus\Z f, \end{equation}
where $S:=\sigma(B)$ and $f$ is the class of an elliptic fiber. Since the fibers are integral, Wazir's Shioda--Tate formula gives $\operatorname{rank}\MW(X/B)=\operatorname{rank}\Pic(X)-\operatorname{rank}\Pic(B)-1$ for the group $\MW(X/B)$ of rational sections \cite[Theorem~3.3]{Wazir}. The first integral hypothesis in \eqref{eq:lattices} therefore forces the global Mordell--Weil rank to be zero.

To obtain sections of positive height, we restrict to a smooth connected curve $D\subset B$ whose preimage $Y:=\pi^{-1}(D)$ is smooth. Set
\begin{equation}\label{eq:notation} p:=\pi|_Y:Y\to D,\quad i:Y\hookrightarrow X,\quad O:=\sigma(D),\quad g:=g(D),\quad \chi_Y:=\chi(\OO_Y)=-K_B\cdot D>0. \end{equation}
For $x\in D$, write $f_x=p^{-1}(x)$ for the scheme-theoretic fiber. Its class is $F$ on $Y$ and $f=i_*F$ on $X$. By a section $C\subset Y$ of $p$, we mean the image of a morphism $D\to Y$ right inverse to $p$. Under fiberwise addition on the generic fiber, the sections form the Mordell--Weil group $\MW(Y/D)$ with identity $O$. We write $\langle\, ,\,\rangle$ for the height pairing on this group and call $\langle C,C\rangle$ the height of $C$ \cite[Definition~8.5]{Shioda}. Since every fiber of $p$ is integral, Shioda's formula \cite[Theorem~8.6]{Shioda} gives
\begin{equation}\label{eq:r} r(C):=\tfrac12\langle C,C\rangle=\chi_Y+C\cdot O, \end{equation}
and in particular $r(O)=0$. By Proposition~\ref{prop:height}, $r(C)$ is the fiber degree of $C$ in $X$. Section~\ref{sec:height} compares this degree with the Jacobi index for a global section and distinguishes it from the discrete charge in Schimannek's GV series.

To record the $n$ added fibers and their multiplicities, let $Y^{[n]}$ denote the Hilbert scheme of length-$n$ zero-dimensional subschemes of $Y$, let $D^{(n)}:=\Sym^nD$ be the $n$-th symmetric product of $D$ (with $Y^{[0]}=D^{(0)}$ a single reduced point), and let $\Pic^d(D)$ be the Picard variety of degree-$d$ line bundles on $D$. Pushing the zero-cycle $[Z]=\sum_{z\in\operatorname{Supp}(Z)}\length(\OO_{Z,z})[z]$ forward along $p$ defines the morphism
\[ \rho_n:Y^{[n]}\longrightarrow D^{(n)},\qquad Z\longmapsto p_*[Z]. \]

We now pass from these cycles to the sheaves they support. For an effective class $\beta\in H_2(X,\Z)$, let $M_\beta$ be Simpson's projective moduli scheme \cite[Theorem~1.21]{Simpson} of Gieseker stable pure one-dimensional sheaves $\mathcal F$ on $X$ with class $\beta$ and $\chi(\mathcal F)=1$. Purity means that every nonzero subsheaf has one-dimensional support, and we write $[\mathcal F]$ for the fundamental cycle weighted by generic lengths. Since $\chi(\mathcal F)=1$, stability is independent of the ample polarization and is equivalent to requiring $\chi(\mathcal F')\leq0$ for every proper nonzero subsheaf $\mathcal F'\subset\mathcal F$. In particular, semistability coincides with stability \cite[Remark~3.2]{MaulikToda}. Sending a sheaf to its support cycle defines the Hilbert--Chow morphism
\[ h_\beta:(M_\beta)_{\mathrm{red}}\longrightarrow\Chow_\beta(X),\qquad \mathcal F\longmapsto[\mathcal F], \]
from the underlying reduced scheme $(M_\beta)_{\mathrm{red}}$ to the Chow variety $\Chow_\beta(X)$ of effective one-cycles of class $\beta$ \cite[Corollary~7.15]{Rydh}\cite[Remark~3.6]{MaulikToda}.

The role of $h_\beta$ in the GV count comes from the BPS interpretation in terms of spins \cite{GV,KKV}. Hosono, Saito, and Takahashi use the relative Lefschetz action for the support morphism to define the invariants \cite{HST}. Kiem and Li introduce vanishing cycles into this construction \cite{KiemLi}. We use Maulik and Toda's definition, which extracts the GV numbers from the direct image of vanishing cycles under $h_\beta$ \cite[Definition~3.7]{MaulikToda}.

To specify the vanishing cycles, fix a trivialization of $K_X$. Then $M_\beta$ carries a natural $d$-critical structure \cite[Corollary~3.19]{BBBJ}\cite[Theorem~3.3]{MaulikToda} recording local presentations as critical loci of regular functions on smooth varieties \cite{Joyce}. An orientation is a line bundle on $(M_\beta)_{\mathrm{red}}$ whose square is the virtual canonical line bundle \cite[Definition~2.31]{Joyce}, and it determines a global perverse sheaf of vanishing cycles $\phi_{M_\beta}$ by Brav, Bussi, Dupont, Joyce, and Szendr\H{o}i \cite[Theorem~6.9]{BBDJS}.

Following \cite[Definition~A.4]{MaulikToda}, an orientation on a proper open and closed subscheme $M\subset M_\beta$ is Calabi--Yau with respect to $h_\beta$ if its line bundle is pulled back from the intermediate variety in the Stein factorization of $h_\beta|_{M_{\mathrm{red}}}$. For such an orientation, the GV contributions $n^M_h$ are defined by
\begin{equation}\label{eq:gv-definition} \sum_{a\in\Z}\chi\bigl(\pH^a Rh_{\beta*}\phi_M\bigr)y^a =\sum_{h\geq0}n^M_h\,(y+2+y^{-1})^h, \end{equation}
where $\pH^a$ denotes middle perverse cohomology and $\chi(\mathcal K):=\sum_{b\in\Z}(-1)^b\dim\HH^b(\Chow_\beta(X),\mathcal K)$. When $M$ is smooth and endowed with its canonical orientation $K_M$, one has $\phi_M=\Q_M[\dim M]$. The following theorem identifies such a component for each infinitesimally rigid section and each $n$.

\begin{theorem}\label{thm:main}
Let $C\subset Y$ be a section of $p:Y\to D$ whose normal bundle $N_{C/X}$ satisfies $H^0(C,N_{C/X})=0$. For every integer $n\geq0$, the locus of stable sheaves with support cycle
\begin{equation}\label{eq:cycle} [C]+\sum_{x\in D}m_x[f_x],\qquad m_x\geq0,\quad\sum_{x\in D}m_x=n, \end{equation}
forms a smooth projective connected component $M_{C,n}$ of $M_{\sigma_*[D]+(r(C)+n)f}$, equipped with an isomorphism
\begin{equation}\label{eq:main-moduli} M_{C,n}\cong\Pic^{g+r(C)}(D)\times Y^{[n]}. \end{equation}
Under this identification, the support morphism factors through the projection onto $Y^{[n]}$ and sends $(N,Z)$ to $[C]+\sum_{x\in D}m_x[f_x]$, where $\sum_{x\in D}m_x[x]=\rho_n(Z)$. Moreover, the canonical orientation on $M_{C,n}$ is Calabi--Yau with respect to the support morphism. Let $n^C_{h,n}$ be the genus-$h$ GV contribution of $M_{C,n}$ in \eqref{eq:gv-definition}, and put $w:=y+2+y^{-1}$. With $q$ recording the number of added fibers, we have
\begin{equation}\label{eq:main-product} \sum_{n,h\geq0}n^C_{h,n}w^hq^n =w^g\prod_{m\geq1} \frac{(1+yq^m)^{2g-2}(1+y^{-1}q^m)^{2g-2}} {(1-q^m)^{12\chi_Y+4g-4}}. \end{equation}
In particular, $n^C_{h,n}=0$ unless $g\leq h\leq g+n$, and
\begin{equation}\label{eq:first-terms} \sum_{h\geq0}n^C_{h,0}w^h=w^g,\qquad \sum_{h\geq0}n^C_{h,1}w^h=12\chi_Yw^g+(2g-2)w^{g+1}. \end{equation}
\end{theorem}

For a fixed pair $(D,C)$, the section starts in fiber degree $r(C)$, so its contribution with $n$ added fibers occurs in degree $r(C)+n$. In \eqref{eq:main-product} the Picard variety contributes $w^g$, while the maps $\rho_n$ record the added fibers in $D^{(n)}$ and give the infinite product. In this sense \eqref{eq:main-product} is the pointwise version of Conjecture~\ref{conj:family}, in which $(D,C)$ varies.

For a rational elliptic surface, the fixed-section product already appears in the relative Lefschetz calculation of Hosono, Saito, and Takahashi \cite[Theorem~4.7]{HST}. In the compact threefold one must also determine when this surface calculation gives a contribution to the ambient GV invariant. Infinitesimal rigidity makes the fixed-section locus a smooth component of the threefold moduli space, on which we determine the orientation. Over a single line in $\PP^2$, Example~\ref{ex:plane-line} gives six rigid sections in the same ambient class. For very general coefficients in this example, Corollary~\ref{cor:rank-two} computes the contribution of every nonzero section in the rank-two subgroup they generate.

There is a related local calculation. When all singular fibers of $Y\to D$ are nodal, Bryan and Kool compute the connected unweighted Donaldson--Thomas series of the total space of $K_Y$ \cite[Corollary~2]{BryanKool}. Setting their variable $p=-y$ gives $(-w)^{g-1}$ times the product in \eqref{eq:main-product}.

Finally, when $D$ is a $(-1)$-curve, every section is infinitesimally rigid and $|D|$ is a point. We can therefore sum the fixed-section contributions to obtain the complete GV series in classes $\sigma_*[D]+kf$. Section~\ref{sec:exceptional} identifies the resulting $E_8$ product with the one-E-string elliptic genus.

Sections~\ref{sec:height} and~\ref{sec:moduli} relate the heights of sections to their ambient curve classes and describe the moduli space over a fixed divisor. Section~\ref{sec:direct} computes the perverse generating series, Section~\ref{sec:obstructions} proves the rigidity and orientation statements, and Section~\ref{sec:exceptional} gives the $E_8$ specialization.

\section{Heights of sections and ambient curve classes}\label{sec:height}

On the smooth surface $Y$, we regard curves as divisors. Let $\NS(Y)$ denote the N\'eron--Severi group of line bundles on $Y$ modulo algebraic equivalence, endowed with the intersection pairing. Relative duality for $p:Y\to D$, adjunction on $B$ and $X$, and the zero section $O=S|_Y$ yield the surface identities
\begin{equation}\label{eq:surface} \begin{gathered} R^1p_*\OO_Y=K_B|_D,\qquad K_Y=N_{Y/X}=p^*N_{D/B},\\ O^2=-\chi_Y,\quad O\cdot F=1,\quad F^2=0,\qquad D^2-\chi_Y=2g-2. \end{gathered} \end{equation}
Since $\deg(K_B|_D)=-\chi_Y<0$, we have $H^0(D,K_B|_D)=0$. The Leray spectral sequence for $p$ therefore gives $h^1(Y,\OO_Y)=g$, and pullback along $p$ induces an isomorphism on identity components of Picard schemes,
\begin{equation}\label{eq:pic} p^*:\Pic^0(D)\xrightarrow{\sim}\Pic^0(Y), \end{equation}
with inverse given by restriction to $O$.

\begin{proposition}\label{prop:height}
For a section $C\in\MW(Y/D)$, define
\[\varphi(C):=[C]-[O]-r(C)F,\qquad \Lambda:=\bigl\{A\in\NS(Y)\;\big|\;A\cdot O=A\cdot F=0\bigr\}.\]
Then $\varphi$ gives an isomorphism of groups $\MW(Y/D)\xrightarrow{\sim}\Lambda$ under which the height pairing takes the form $\langle C_1,C_2\rangle=-\varphi(C_1)\cdot\varphi(C_2)$. Moreover,
\begin{equation}\label{eq:height-class} \NS(Y)=\langle O,F\rangle\oplus\Lambda,\qquad \varphi(C)^2=-2r(C),\qquad i_*[C]=\sigma_*[D]+r(C)f. \end{equation}
The group $\MW(Y/D)$ is torsion-free and contains only finitely many sections of bounded height, and every nonzero section $C\ne O$ satisfies $r(C)\geq\chi_Y$.
\end{proposition}
\begin{proof}
Since every fiber of $p$ is integral, the local correction terms from reducible fibers in Shioda's height formula vanish, giving $\langle C,C\rangle=2\chi_Y+2C\cdot O$ \cite[Theorem~8.6]{Shioda}. For $C\ne O$, we have $C\cdot O\geq0$ and hence $r(C)\geq\chi_Y>0$, which shows that $\MW(Y/D)$ has no torsion. The Shioda--Tate theorem identifies $\NS(Y)/\langle O,F\rangle$ with $\MW(Y/D)$ \cite[Theorem~1.3]{Shioda}.

The classes $O,F$ span a free subgroup, and the quotient $\MW(Y/D)$ is torsion-free. Hence $\NS(Y)$ is torsion-free as well.

With respect to the basis $O,F$, the intersection form on $\langle O,F\rangle$ has matrix
\[\begin{pmatrix}-\chi_Y&1\\1&0\end{pmatrix},\]
which has determinant $-1$ and signature $(1,1)$. Orthogonal projection off $\langle O,F\rangle$ is therefore defined over $\Z$, yielding the direct sum decomposition $\NS(Y)=\langle O,F\rangle\oplus\Lambda$ together with the projection $\varphi$ and self-intersection $\varphi(C)^2=-2r(C)$. By the Hodge index theorem, the intersection form on $\Lambda$ is negative definite, so only finitely many elements of $\MW(Y/D)$ can have bounded height.

Finally, for any class $A\in\NS(Y)$, combining the projection formula with \eqref{eq:lattices} gives
\begin{equation}\label{eq:pushforward} i_*A=(A\cdot F)\sigma_*[D] +\bigl(A\cdot O+\chi_Y(A\cdot F)\bigr)f. \end{equation}
Indeed, pairing $i_*A$ with pullbacks of divisors from $B$ determines the coefficient of $\sigma_*[D]$, while intersecting with $S$ and using $S\vert{}_S=K_B$ determines the coefficient of $f$. It follows that $i_*\Lambda=0$, and evaluating \eqref{eq:pushforward} on $[C]$ yields $i_*[C]=\sigma_*[D]+r(C)f$.
\end{proof}

Proposition~\ref{prop:height} applies to a section of $Y$ defined only over $D$. To compare with the height divisor of a global section, let $\widetilde\pi:\widetilde X\to B$ be a smooth elliptic Calabi--Yau threefold with zero section $S_0$ and no reducible fibers in codimension one. Wazir defines a pairing of divisors on $\widetilde X$ with values in $\Pic(B)$ by $(A,A')\mapsto\widetilde\pi_*(A\cdot A')$ \cite[Definition~3.2]{Wazir}. Write $P$ also for the divisor of a section. In the convention of Morrison and Park \cite[equations~(3.4)--(3.8)]{MP}, the Shioda class of $P$ and the height pairing of $P,Q$ are
\[
\Sigma(P)=P-S_0-\widetilde\pi^*\widetilde\pi_*((P-S_0)\cdot S_0),\qquad b(P,Q)=-\widetilde\pi_*(\Sigma(P)\cdot\Sigma(Q)).
\]
Thus the divisor-valued height pairing is the negative of Wazir's pairing on the Shioda classes. Put $b(P)=b(P,P)$. If $\widetilde Y=\widetilde\pi^{-1}(D)$ is smooth with integral fibers and $P_D=P|_D$, then $\Sigma(P)|_{\widetilde Y}=\varphi(P_D)$ in $\NS(\widetilde Y)$. The projection formula gives
\[
\tfrac12 D\cdot b(P)=-\tfrac12\bigl(\Sigma(P)|_{\widetilde Y}\bigr)^2=r(P_D).
\]
The number $D\cdot b(P)/2$ is the Jacobi index of the charge variable in base class $[D]$ \cite[equation~(3.45)]{LLW}. Fierro Cota, Klemm, and Schimannek give the corresponding index matrix for several sections \cite[equations~(4.32)--(4.33)]{FCKS}. Schimannek uses the convention $b=-b(P)$ and writes the same index as $-D\cdot b/2$ \cite[equations~(4.3), (4.8)]{SchDiscrete}.

In \eqref{eq:lattices}, $H_2(X,\Z)$ is torsion-free and $r(C)$ records the ordinary fiber degree. The discrete charges in Schimannek's genus-one examples arise from torsion curve classes, weighted by characters of a flat $B$-field \cite[equations~(3.12)--(3.13)]{SchDiscrete}. The Tate--Shafarevich group parametrizes the locally trivial genus-one torsors of the Jacobian \cite{DG}. In Schimannek's example over $\PP^2$, the two $\Z/2$ sectors in fiber degree one and base degree zero sum to the smooth-model invariant $284+256=540$ \cite[\S6.3, Tables~6--7]{SchDiscrete}. Toda computes the local GV invariants in every pure fiber class of the smooth Weierstrass model \cite[Theorem~3.5]{TodaWall}. For the torsion grading, Katz, Klemm, Schimannek, and Sharpe propose a definition using stable sheaves on small resolutions of the nodal model \cite[\S3.5]{KKSS}.

\section{Moduli of stable sheaves over a fixed curve}\label{sec:moduli}

We use the relative Fourier--Mukai equivalence of Bridgeland and Maciocia \cite{BM}, with the normalization and moduli comparison of Diaconescu \cite[\S3.1 and Theorem~1.1]{Diaconescu}. Let $p_1,p_2:X\times_BX\to X$ be the two projections and $\Delta\subset X\times_BX$ the relative diagonal with ideal sheaf $I_\Delta$. Define the normalized kernel and integral transform by
\begin{equation}\label{eq:FM} \mathcal P:=I_\Delta\otimes p_1^*\OO_X(S)\otimes p_2^*\OO_X(S)\otimes (\pi\circ p_1)^*K_B^{-1},\qquad \Phi(-):=Rp_{2*}\bigl(Lp_1^*(-)\otimes^{\mathbf L}\mathcal P\bigr). \end{equation}
Here $X$ is identified with the relative compactified degree-zero Picard scheme over $B$, which parametrizes rank-one torsion-free sheaves of degree zero on the fibers, via $x\mapsto I_x\otimes\OO_{X_{\pi(x)}}(S)$, where $I_x$ is the ideal of $x$ in its fiber. On a smooth fiber over $b\in B$, this is $x\mapsto\OO_{X_b}(\sigma(b)-x)$. See \cite[\S3.1]{Diaconescu} for this normalization and the resulting derived equivalence.

For $k\geq0$, put $\beta_k:=\sigma_*[D]+kf$. Diaconescu's comparison theorem \cite[Theorem~1.1 and equations~(3.8)--(3.9)]{Diaconescu} identifies $M_{\beta_k}$ with the moduli space of adiabatically semistable pure two-dimensional sheaves $E$ on $X$ with Chern character
\begin{equation}\label{eq:characters} \ch_1(E)=[Y],\qquad \ch_2(E)=(1-\chi_Y/2)f,\qquad \ch_3(E)=-k, \end{equation}
where we identify $\ch_1$ and $\ch_2$ with homology classes via Poincar\'e duality and view $\ch_3$ as a rational number by taking degree. Concretely, if we fix polarizations $\widehat\omega=S-s\pi^*K_B$ and $\omega_t=tS-s\pi^*K_B$ with $s\gg0$ and $0<t\ll1$ chosen relative to the numerical class, adiabatic semistability means Gieseker semistability with respect to $\omega_{t'}$ for all $0<t'\leq t$. This comparison is an isomorphism of moduli stacks compatible with flat families, while $M_{\beta_k}$ itself does not depend on $\widehat\omega$. The pushforward $i_*J$ of any rank-one torsion-free sheaf $J$ on the smooth integral surface $Y$ is Gieseker stable with respect to every ample divisor on $X$, and hence belongs to the adiabatically stable locus whenever its Chern character is \eqref{eq:characters}.

Next, we construct the projection to the linear system $|D|:=\PP H^0(B,\OO_B(D))$,
\begin{equation}\label{eq:divisor-map} s_k:M_{\beta_k}\longrightarrow|D|, \qquad M_k(D):=s_k^{-1}(D), \end{equation}
where $M_k(D)$ denotes the scheme-theoretic fiber over $D\in|D|$. For a family $\mathcal F_T$ of stable sheaves on $X\times T$ flat over $T$, the relative transform with kernel \eqref{eq:FM} gives a $T$-flat family of pure codimension-one sheaves on $X$. The divisorial support of each fiber is the sum of its support divisors weighted by the lengths at their generic points. In families these divisors form a relative effective Cartier divisor on $X\times T$ whose formation commutes with base change \cite[\S11, Corollary~11.8]{Rydh}\cite[Lemmas~7.29 and~7.32]{Kollar}. Since $H^1(X,\OO_X)=0$, the line bundle associated with this relative divisor agrees with $\pi^*\OO_B(D)$ up to twist by a line bundle from $T$. Because $H^0(X,\pi^*\OO_B(D))=H^0(B,\OO_B(D))$, the divisor is pulled back from a unique relative effective divisor on $B\times T$ in $|D|$. The construction descends to $s_k$ because it commutes with base change and twisting a family by a line bundle from $T$. To identify its value at $\mathcal F$, localize at the generic point of each component $\Gamma$ of $\pi_*[\mathcal F]$. Filter the localized finite-length sheaf by powers of a uniformizer of $\OO_{B,\eta_\Gamma}$. Over an algebraic closure of the residue field, the graded pieces have a total of $\operatorname{mult}_\Gamma(\pi_*[\mathcal F])$ point factors. Each point factor transforms to a rank-one sheaf on the elliptic fiber \cite[Lemma~3.1]{Diaconescu}. Hence the multiplicity of $\pi^{-1}(\Gamma)$ in the divisorial support is the same. Components over closed points of $B$ contribute no divisor, so $s_k(\mathcal F)=\pi_*[\mathcal F]$.

On an elliptic surface with integral fibers, Yoshioka uses the shift $d(C)=(C-O)\cdot O=r(C)$ and a subscheme of length $k-d(C)$ in his Fourier--Mukai correspondence \cite[\S3.2, equations~(3.13)--(3.15) and Theorem~3.15]{Yoshioka}. For a rational elliptic surface, Hosono, Saito, and Takahashi give the support morphism in a fixed section class and attribute the description to Yoshioka \cite[Proposition~4.8]{HST}. We identify the scheme-theoretic fiber in the compact threefold moduli space.

\begin{theorem}\label{thm:fixed}
The scheme-theoretic fiber $M_k(D)$ in \eqref{eq:divisor-map} is smooth and decomposes as
\begin{equation}\label{eq:fixed-product} M_k(D)\cong\coprod_{\substack{C\in\MW(Y/D)\\r(C)\leq k}} \Pic^{g+r(C)}(D)\times Y^{[k-r(C)]}. \end{equation}
Under this isomorphism, a point $(N,Z)$ in the component indexed by $C$ corresponds to a stable sheaf $\mathcal F$ with Fourier--Mukai transform
\begin{equation}\label{eq:rankone} \Phi(\mathcal F)=i_*J,\qquad J=I_Z\otimes\OO_Y(O-C)\otimes p^*N, \quad\deg N=g+r(C),\quad\length Z=k-r(C), \end{equation}
where $I_Z\subset\OO_Y$ is the ideal sheaf of $Z$. The multiplicities in the support cycle $[\mathcal F]$ of \eqref{eq:cycle} are determined by $\sum_{x\in D}m_x[x]=\rho_{k-r(C)}(Z)=p_*[Z]$.
\end{theorem}
\begin{proof}
Let $\mathcal F$ be a geometric point of $M_k(D)$ and set $E:=\Phi(\mathcal F)$. Because the divisorial support of $E$ is $Y$ with multiplicity one, $E$ has length one at the generic point of $Y$ and is annihilated there by any local defining equation of $Y$. Multiplication on $E$ by a local equation of $Y$ therefore has image supported in dimension at most one. Purity forces this image to vanish, so $E=i_*J$ for a rank-one torsion-free sheaf $J$ on $Y$. Since $Y$ is a smooth surface, $L:=J^{\vee\vee}$ is a line bundle and $J\cong I_Z\otimes L$ for a unique zero-dimensional subscheme $Z\subset Y$.

Put $\delta:=c_1(L)$, $\nu:=c_1(N_{Y/X})=D^2F$, and $n:=\length Z$. Applying Grothendieck--Riemann--Roch to $i:Y\hookrightarrow X$ yields
\begin{equation}\label{eq:GRR} \begin{aligned} \ch_2(i_*J)&=i_*(\delta-\nu/2),\\ \ch_3(i_*J)&=\int_Y\left(\delta^2/2-n-\delta\nu/2+\nu^2/6\right). \end{aligned} \end{equation}
Combining \eqref{eq:characters}, \eqref{eq:surface}, and \eqref{eq:pushforward}, we find
\[i_*\delta=\left(1-\frac{\chi_Y}{2}+\frac{D^2}{2}\right)f=gf, \qquad \delta\cdot F=0.\]
By Proposition~\ref{prop:height}, there exists a unique section $C\in\MW(Y/D)$ such that
\[\delta=-\varphi(C)+gF=O-C+(g+r(C))F,\]
which by \eqref{eq:pic} is equivalent to $L\cong\OO_Y(O-C)\otimes p^*N$ for a unique line bundle $N\in\Pic^{g+r(C)}(D)$. Furthermore,
\[\delta^2=-2r(C),\qquad\delta\nu=\nu^2=0,\qquad -k=-r(C)-n,\]
which establishes the classification of closed points together with the numerical relations in \eqref{eq:rankone}. Conversely, every sheaf of the form \eqref{eq:rankone} has Chern character \eqref{eq:characters} and is Gieseker stable for every polarization, so its inverse transform defines a point of $M_k(D)$.

Now choose a Poincar\'e line bundle on $D\times\Pic^{g+r(C)}(D)$ and the universal ideal sheaf on $Y\times Y^{[n]}$. The resulting family gives a morphism to $M_{\beta_k}$. Since the source is reduced and $s_k$ maps every geometric point to $D$, this morphism factors through $M_k(D)$. The map from the disjoint union of the products is proper and bijective on geometric points. To check the scheme structure, consider the change of rings sequence
\begin{equation}\label{eq:normal-ext} \begin{aligned} 0\longrightarrow\Ext^1_Y(J,J)&\longrightarrow\Ext^1_X(i_*J,i_*J) \xrightarrow{\eta}H^0(Y,N_{Y/X})\\ &\xrightarrow{o_J}\Ext^2_Y(J,J), \end{aligned} \end{equation}
where $\mathcal{H}om_Y(J,J)=\OO_Y$ identifies the normal term. Since $N_{Y/X}=p^*N_{D/B}$ and $H^1(B,\OO_B)=0$,
\[ H^0(Y,N_{Y/X})=H^0(D,N_{D/B})\cong T_D|D|. \]
On $Y\setminus Z$, the sheaf $J$ is invertible. In a first-order deformation, a local equation $z$ of $Y$ acts on a lift of a local generator as multiplication by $\varepsilon a$, where $\varepsilon^2=0$. The normal map $\eta$ sends this deformation to $a$, and its divisorial support has equation $z-\varepsilon a$. Hence $\eta$ and $ds_k$ agree on $Y\setminus Z$. They agree globally as sections of $N_{Y/X}$, giving $T_{\mathcal F}M_k(D)=\ker\eta=\Ext^1_Y(J,J)$.

Since $\mathcal{H}om_Y(J,J)=\OO_Y$, Serre duality gives $\dim\Ext^2_Y(J,J)=h^2(\OO_Y)$. Riemann--Roch gives $\sum_{j=0}^2(-1)^j\dim\Ext^j_Y(J,J)=\chi_Y-2n$. With $\dim\Hom_Y(J,J)=1$, we obtain
\[\dim\Ext^1_Y(J,J)=g+2n.\]
By Fogarty's theorem \cite[Theorem~2.4]{Fogarty}, each product $\Pic^{g+r(C)}(D)\times Y^{[n]}$ is smooth of dimension $g+2n$. The family map is finite by properness and quasi-finiteness. It is a homeomorphism, so the local dimension $g+2n$ at each closed point equals the tangent dimension. Thus $M_k(D)$ is regular. In characteristic zero the map is birational on each component, hence an isomorphism by normality.

It remains to determine the support cycle. Write $\Psi=\Phi^{-1}$. Restricting the kernel \eqref{eq:FM} along $C$ gives $\OO_Y(O-C)$ up to a pullback from $D$. Thus $\Psi(i_*L)$ is a line bundle on $C$, extended by zero to $X$. For $z\in Y$, let $k_z$ be its skyscraper sheaf. Then $\Psi(k_z)[-1]$ is rank one and torsion-free on $f_{p(z)}$ \cite[Lemma~3.1 and equation~(3.5)]{Diaconescu}. A composition series of $L|_Z$ shows that $\Psi(i_*(L|_Z))[-1]$ is a sheaf filtered by the $\Psi(k_z)[-1]$, each with multiplicity $\length(\OO_{Z,z})$. Applying $\Psi$ to $0\to i_*J\to i_*L\to i_*(L|_Z)\to0$ therefore gives a short exact sequence
\[ 0\longrightarrow\Psi(i_*(L|_Z))[-1]\longrightarrow\mathcal F\longrightarrow\Psi(i_*L)\longrightarrow0. \]
Additivity of fundamental cycles gives $[\mathcal F]=[C]+\sum_{z\in\operatorname{Supp}(Z)}\length(\OO_{Z,z})[f_{p(z)}]$, as claimed.
\end{proof}

\section{Perverse direct images and generating series}\label{sec:direct}

Recall that $\rho_n:Y^{[n]}\to D^{(n)}$ sends a length-$n$ subscheme $Z$ to the divisor $p_*[Z]$ under the elliptic fibration $p:Y\to D$. Define
\[\mathcal A_n:=R\rho_{n*}\Q_{Y^{[n]}}[2n],\qquad Q_n(y):=\sum_{a\in\Z}\chi\bigl(\pH^a\mathcal A_n\bigr)y^a, \qquad\mathcal A_0:=\Q,\quad Q_0(y):=1.\]
Thus $Q_n(y)$ records the Euler characteristics on $D^{(n)}$ of the perverse cohomology sheaves of $\mathcal A_n$. Hosono, Saito, and Takahashi discuss the extension of their rational elliptic surface calculation to other elliptic surfaces in the remark following the proof of \cite[Theorem~4.7]{HST}.

\begin{lemma}\label{lem:hilb-product}
With $g=g(D)$ and $\chi_Y=\chi(\OO_Y)$, one has
\begin{equation}\label{eq:hilb-product} \sum_{n\geq0}Q_n(y)q^n =\prod_{m\geq1} \frac{(1+yq^m)^{2g-2}(1+y^{-1}q^m)^{2g-2}} {(1-q^m)^{12\chi_Y+4g-4}}. \end{equation}
\end{lemma}
\begin{proof}
Let $j:U\hookrightarrow D$ denote the open locus over which $p$ is smooth, and set $\mathbb V:=(R^1p_*\Q_Y)|_U$ with intersection complex $\IC_D(\mathbb V):=j_{!*}\mathbb V[1]$. The decomposition theorem of Beilinson, Bernstein, Deligne, and Gabber \cite[Theorem~6.2.5]{BBD}, applied to $p:Y\to D$, gives
\begin{equation}\label{eq:A1} \mathcal A_1=Rp_*\Q_Y[2] \cong\Q_D[2]\oplus\IC_D(\mathbb V)\oplus\Q_D. \end{equation}
At $x\in D$, the multiplicity of a skyscraper summand is one less than the number of irreducible components of $f_x$, hence zero. This also follows by comparing degree-two stalks, since $H^2(f_x,\Q)=\Q$, and applying Verdier duality.

By the canonical bundle formula for elliptic surfaces and Noether's formula, the topological Euler characteristic of $Y$ is $e(Y)=12\chi_Y$. Since $b_1(Y)=\dim H^1(Y,\Q)=2g$, we obtain $b_2(Y)=12\chi_Y+4g-2$. Taking hypercohomology in \eqref{eq:A1} shows that $\HH^b(D,\IC_D(\mathbb V))$ vanishes for $b\ne0$ and has dimension $12\chi_Y+4g-4$ for $b=0$. The nonzero bigraded hypercohomology spaces $\HH^b(D,\pH^a\mathcal A_1)$ are therefore recorded by
\begin{equation}\label{eq:bigrading} \begin{array}{c|c|c} (a,b)&\dim\HH^b(D,\pH^a\mathcal A_1)&a+b\pmod2\\ \hline (1,1),(1,-1),(-1,1),(-1,-1)&1\text{ for each pair}&0\\ (1,0),(-1,0)&2g\text{ for each pair}&1\\ (0,0)&12\chi_Y+4g-4&0. \end{array} \end{equation}

For a partition $\lambda=(1^{m_1}2^{m_2}\cdots)$ of $n$, consider the finite map
\[ \mu_\lambda:\prod_{j:\,m_j>0}D^{(m_j)}\longrightarrow D^{(n)},\qquad (\xi_j)_j\longmapsto\sum_j j\,\xi_j. \]
Pushing forward the Hilbert--Chow decomposition \cite[Theorem~3 and Proposition~1]{GS} (in the form of \cite[Theorem~4.1.1]{dCM}) to $D^{(n)}$ gives an isomorphism of rational constructible complexes
\begin{equation}\label{eq:partitions} \mathcal A_n\cong\bigoplus_{\lambda\vdash n}R\mu_{\lambda*} \left(\mathop{\boxtimes}_{j\geq1}\Sym^{m_j}\mathcal A_1\right), \end{equation}
where $\Sym^m\mathcal K$ denotes the invariant summand of the finite pushforward of $\mathcal K^{\boxtimes m}$ to $D^{(m)}$ with the usual Koszul signs. We omit Tate twists for rational complexes. For the induced perverse filtration on the Hilbert scheme, see \cite[Proposition~4.12]{Zhang} and \cite[Proposition~2.2]{ShenZhang}.

Since finite pushforwards are perverse exact \cite[\S2.2.6]{BBD}, external products and the K\"unneth formula reduce the calculation to the bigraded vector space in \eqref{eq:bigrading}, with permutation parity $a+b$ and Euler characteristic sign $(-1)^b$. For each part size $m$, the three rows in \eqref{eq:bigrading} contribute $(1+yq^m)^{-2}(1+y^{-1}q^m)^{-2}$, $(1+yq^m)^{2g}(1+y^{-1}q^m)^{2g}$, and $(1-q^m)^{-(12\chi_Y+4g-4)}$, respectively, and multiplying over $m\geq1$ proves \eqref{eq:hilb-product}.
\end{proof}

Setting $y=-1$ in Lemma~\ref{lem:hilb-product} recovers the Euler characteristic specialization of G\"ottsche's formula \cite{Gottsche},
\[ Q_n(-1)=e(Y^{[n]}),\qquad \sum_{n\geq0}e(Y^{[n]})q^n=\prod_{m\geq1}(1-q^m)^{-12\chi_Y}. \]

For a fixed section $C$ and $n$ added fibers, the total fiber degree is $k=r(C)+n$. The morphism
\[ \tau_{C,n}:D^{(n)}\longrightarrow\Chow_{\beta_k}(X),\qquad \sum_xm_x[x]\longmapsto[C]+\sum_xm_x[f_x] \]
is injective on geometric points since the one-cycle recovers its fiber multiplicities, and it is finite by properness and quasi-finiteness. On the component $\Pic^{g+r(C)}(D)\times Y^{[n]}$ of Theorem~\ref{thm:fixed}, the support map is $h^D_{C,n}=\tau_{C,n}\circ\rho_n\circ\operatorname{pr}_2$, where $\operatorname{pr}_2$ is projection to $Y^{[n]}$. The K\"unneth formula gives
\begin{equation}\label{eq:fixed-direct} R(h^D_{C,n})_*\IC \cong (\tau_{C,n})_*\bigl(R\Gamma(\Pic^{g+r(C)}(D),\Q)[g]\otimes\mathcal A_n\bigr). \end{equation}
Here $\IC$ is the constant sheaf on the product shifted by $g+2n$. Since finite pushforward is perverse exact and
\[ \sum_{i=0}^{2g}\dim H^i(\Pic^d(D),\Q)y^{i-g} =y^{-g}(1+y)^{2g}=(y+2+y^{-1})^g=w^g, \]
the contribution of the intersection complex on $\Pic^{g+r(C)}(D)\times Y^{[n]}$ is $w^gQ_n(y)$.

\begin{proposition}\label{prop:fixed-series}
Fix $\beta_k=\sigma_*[D]+kf$. On the fiber $M_k(D)=s_k^{-1}(D)$ over the projected curve $D$, let $\IC_{M_k(D)}$ be the sum of the constant sheaves on the components of \eqref{eq:fixed-product}, each shifted by the dimension of that component. Define
\[ P_{D,k}(y):=\sum_a\chi\bigl(\pH^a R(h_{\beta_k}|_{M_k(D)})_*\IC_{M_k(D)}\bigr)y^a, \qquad \Theta_Y(q):=\sum_{C\in\MW(Y/D)}q^{r(C)}. \]
Then
\begin{equation}\label{eq:fixed-series} \sum_{k\geq0}P_{D,k}(y)q^k =w^g\Theta_Y(q)\prod_{m\geq1} \frac{(1+yq^m)^{2g-2}(1+y^{-1}q^m)^{2g-2}} {(1-q^m)^{12\chi_Y+4g-4}}. \end{equation}
\end{proposition}
\begin{proof}
By Proposition~\ref{prop:height}, the height bound ensures that each coefficient of $\Theta_Y(q)$ is finite. Summing \eqref{eq:fixed-direct} over the components of \eqref{eq:fixed-product} gives
\[ P_{D,k}(y)=w^g\sum_{\substack{C\in\MW(Y/D)\\r(C)\leq k}} Q_{k-r(C)}(y). \]
Multiplying by $q^k$ and summing over $k$ gives \eqref{eq:fixed-series} by Lemma~\ref{lem:hilb-product}.
\end{proof}

Formula \eqref{eq:fixed-series} concerns the intrinsic intersection complexes on the fixed scheme-theoretic fiber $M_k(D)$. When $D$ or $C$ deforms in $X$, this fiber need not be open in $M_{\beta_k}$, and the next section determines when the component calculation computes the ambient GV invariants \eqref{eq:gv-definition}.

\section{Normal deformations and vanishing cycles}\label{sec:obstructions}

\begin{proposition}\label{prop:tangent}
Let $J=I_Z\otimes\OO_Y(O-C)\otimes p^*N$ be the sheaf on $Y$ from \eqref{eq:rankone}, with $n=\length Z$. Then
\begin{equation}\label{eq:tangent} \dim\Ext^1_X(i_*J,i_*J)=g+2n+h^0(C,N_{C/X}). \end{equation}
\end{proposition}
\begin{proof}
The differential of $p$ splits along $C$, so $p$ is smooth there. By adjunction, $N_{C/Y}=K_C\otimes K_Y^{-1}|_C=K_D\otimes N_{D/B}^{-1}=K_B|_D$. The normal bundle sequence for $C\subset Y\subset X$ is therefore
\begin{equation}\label{eq:section-normal} 0\longrightarrow K_B|_D\longrightarrow N_{C/X} \longrightarrow N_{D/B}\longrightarrow0. \end{equation}
Let $\delta_C:H^0(D,N_{D/B})\to H^1(D,K_B|_D)$ be the connecting map. As $H^0(D,K_B|_D)=0$, its kernel is $H^0(C,N_{C/X})$. Moreover, the trace map $\tr:\Ext^2_Y(J,J)\to H^2(Y,\OO_Y)$ is an isomorphism, since its Serre dual is the isomorphism $H^0(Y,K_Y)\to\Hom_Y(J,J\otimes K_Y)$ given by the identity on $J$.

For $v\in H^0(Y,N_{Y/X})$, let $Y_v$ be the corresponding deformation over $\C[\varepsilon]/(\varepsilon^2)$. Choose a local equation $z=0$ for $Y$ and a lift $a$ of $v(z)$, so that $Y_v$ has equation $z-\varepsilon a=0$. For a flat deformation $E$ of $i_*J$, multiplication by $\varepsilon$ identifies $E/\varepsilon E$ with $\varepsilon E$. Under this identification, multiplication by $z$ induces the endomorphism of $J$ given by $\eta(E)(z)$, where $\eta$ is the normal map in \eqref{eq:normal-ext}. Thus $\eta(E)=v$ precisely when $(z-\varepsilon a)E=0$ on each such chart. Equivalently, $E$ is the pushforward of a flat lift of $J$ to $Y_v$. Exactness of \eqref{eq:normal-ext} gives
\[ o_J(v)=0\quad\Longleftrightarrow\quad J\text{ lifts to }Y_v. \]

Write $\kappa(v)\in\Ext^1_Y(\Omega_Y,\OO_Y)$ for the Kodaira--Spencer class of $Y_v$ and $\operatorname{At}(J)\in\Ext^1_Y(J,J\otimes\Omega_Y)$ for the Atiyah class of $J$. Since $Y$ is smooth, $J$ is perfect for every $Z$. Its lifting obstruction is $(1_J\otimes\kappa(v))\circ\operatorname{At}(J)$, whose trace is the obstruction $\kappa(v)\mathbin{\lrcorner}c_1(\det J)$ to lifting $\det J$ \cite[Corollary~3.4 and \S4.2]{HT}. The trace is an isomorphism, so $J$ lifts if and only if $\det J$ lifts. Now
\[ \det I_Z=\OO_Y,\qquad \det J=\OO_Y(O-C)\otimes p^*N. \]
The factor $\OO_Y(O)$ extends from $\OO_X(S)$, and $N$ extends to the induced deformation of $D$ since $H^2(D,\OO_D)=0$. Hence $o_J(v)=0$ precisely when $\OO_Y(C)$ lifts to $Y_v$. We show that this is equivalent to lifting the divisor $C$. Let $s_C$ be its canonical section. On every integral genus-one fiber $\OO_Y(C)$ has degree one, so $s_C$ spans $H^0$ and $H^1=0$. Cohomology and base change give
\[ p_*\OO_Y(C)=\OO_D,\qquad R^1p_*\OO_Y(C)=0,\qquad H^1(Y,\OO_Y)\xrightarrow{\,s_C\,}H^1(Y,\OO_Y(C))\ \text{is an isomorphism}. \]
Indeed, $H^0(D,R^1p_*\OO_Y)=0$ by \eqref{eq:surface}.

If $\OO_Y(C)$ lifts to $Y_v$ over $\C[\varepsilon]/(\varepsilon^2)$, choose a lift $L_v$. The obstruction to lifting $s_C$ lies in $H^1(Y,\OO_Y(C))$. Twisting $L_v$ by a line bundle with transition functions $1+\varepsilon u_{ij}$ changes that obstruction by the class of $(u_{ij}s_C)$. The displayed isomorphism allows us to cancel it. The lifted section cuts out a relative Cartier divisor lifting $C$. Conversely, a lifted Cartier divisor gives a lift of $\OO_Y(C)$.

The normal sequence \eqref{eq:section-normal} says that $C$ lifts inside $Y_v$ precisely when $\delta_C(v)=0$. Hence
\begin{equation}\label{eq:obstruction} \ker o_J=\ker\delta_C\cong H^0(C,N_{C/X}). \end{equation}
Substitution into \eqref{eq:normal-ext} proves \eqref{eq:tangent}.
\end{proof}

\begin{proof}[Proof of Theorem~\ref{thm:main}]
By Theorem~\ref{thm:fixed}, the product $M_{C,n}$ in \eqref{eq:main-moduli} is a smooth projective connected component of the moduli fiber $M_{r(C)+n}(D)$. Since the Fourier--Mukai equivalence preserves $\Ext^1$, Proposition~\ref{prop:tangent} and the vanishing $H^0(C,N_{C/X})=0$ imply that the ambient tangent dimension at each closed point of $M_{C,n}$ is $g+2n$. At such a point the local dimension lies between $\dim M_{C,n}=g+2n$ and the tangent dimension. The ambient local ring is therefore regular. The ideal of $M_{C,n}$ in each such regular local ring is zero, since a nonzero ideal in a local domain lowers the dimension of the quotient. Thus the closed immersion $M_{C,n}\hookrightarrow M_{[C]+nf}$ is also open and identifies $M_{C,n}$ with a connected component of $M_{[C]+nf}$.

Let $\operatorname{pr}_2:M_{C,n}\to Y^{[n]}$ be the second projection. The line bundle $N_{D/B}^{\boxtimes n}$ on $D^n$ descends to a line bundle $N_{D/B}^{(n)}$ on $D^{(n)}$ because point stabilizers act trivially on its fibers. Since the Picard variety has trivial canonical bundle, $K_Y=p^*N_{D/B}$, and the Hilbert--Chow morphism is crepant, the local construction of the canonical bundle of $Y^{[n]}$ in \cite[\S6]{Beauville} shows that $K_{Y^{[n]}}$ is the pullback of the symmetric descent of $K_Y^{\boxtimes n}$ on $Y^{(n)}$, and hence
\begin{equation}\label{eq:canonical} K_{M_{C,n}}=\operatorname{pr}_2^*K_{Y^{[n]}} =(\rho_n\circ\operatorname{pr}_2)^*N_{D/B}^{(n)}. \end{equation}

The map $\rho_n\circ\operatorname{pr}_2$ sends $(N,Z)$ to $p_*[Z]$. Over $x_1+\cdots+x_n$ with distinct $x_i\in D$ and smooth fibers $p^{-1}(x_i)$, its fiber is $\Pic^{g+r(C)}(D)\times\prod_i p^{-1}(x_i)$. Thus the proper surjective map has geometrically connected generic fiber, also for $n=0$. Since its source is smooth and connected and its target $D^{(n)}$ is normal, Stein factorization gives
\[ (\rho_n\circ\operatorname{pr}_2)_*\OO_{M_{C,n}}=\OO_{D^{(n)}}. \]
The map $\tau_{C,n}$ from $D^{(n)}$ to the Chow variety adds $[C]$ to the fiber cycle. It is finite, so this is the Stein factorization of the support morphism. On a smooth $d$-critical scheme, the virtual canonical line bundle is $K_{M_{C,n}}^{\otimes2}$ \cite[Theorem~2.28]{Joyce}, and \eqref{eq:canonical} shows that its square root $K_{M_{C,n}}$ descends to the Stein factor and therefore defines a Calabi--Yau orientation in the sense of \cite[Definition~A.4]{MaulikToda}. With this orientation, the smooth critical chart with zero potential gives \cite[Theorem~6.9]{BBDJS}
\[ \phi_{M_{C,n}}=\Q_{M_{C,n}}[g+2n]. \]
Consequently, the left-hand side of \eqref{eq:gv-definition} is $w^gQ_n(y)$, and Lemma~\ref{lem:hilb-product} proves \eqref{eq:main-product}.

Finally, each $Q_n(y)$ is a symmetric Laurent polynomial with exponents between $-n$ and $n$, and hence is an integral polynomial in $w=y+2+y^{-1}$ of degree at most $n$. Multiplying by $w^g$ gives the genus bound, while extracting the coefficient of $q$ in \eqref{eq:hilb-product} gives $12\chi_Y+(2g-2)w$, proving \eqref{eq:first-terms}.
\end{proof}

\begin{remark}\label{rem:plane}
For $B=\PP^2$ and a smooth plane curve $D$ of degree $d$, we have $g=(d-1)(d-2)/2$ and $\chi_Y=3d$. In this case, the vanishing condition in Theorem~\ref{thm:main} is equivalent to the injectivity of
\[ H^0(D,\OO_D(d))\xrightarrow{\delta_C}H^1(D,\OO_D(-3)) \]
from \eqref{eq:section-normal}, where both source and target have dimension $d(d+3)/2$. By Proposition~\ref{prop:height}, every nonzero section satisfies $r(C)\geq3d$, which for $d=1$ matches the threshold for additional sections in HKKW \cite[\S5.1.7]{HKKW}.
\end{remark}

\begin{example}\label{ex:plane-line}
Let $[s:t:u]$ be homogeneous coordinates on $B=\PP^2$ and let $D=(u=0)$. Choose Weierstrass coefficients
\[ a=t^{12}+uA_{11},\qquad b=s^{18}+us^7t^{10}+u^2B_{16}, \]
where $A_{11}$ and $B_{16}$ are general homogeneous polynomials of the indicated degrees. Along $D$, the discriminant restricts to $4t^{36}+27s^{36}$, which has $36$ simple zeros, and hence $D$ meets the discriminant transversely. Off $D$, varying $A_{11}$ and $B_{16}$ prescribes the first jets of $a,b$ independently. For general coefficients the pair $(a,b)$ is transverse to the smooth part of $4a^3+27b^2=0$ and to its cusp point. Thus $X$ is smooth with the required discriminant, and $Y=\pi^{-1}(D)$ is smooth.

The smooth Weierstrass models form an open subset of the affine space of coefficients. Its general member is the crepant resolution of a degree-$18$ hypersurface in $\PP(1,1,1,6,9)$ \cite[\S2]{HKK2021}. The five rays corresponding to the vertices generate the toric lattice $\Z^5/\Z(1,1,1,6,9)$. The formulas of Batyrev and Kreuzer \cite[Theorems~1.6 and~3.8]{BatyrevKreuzer} therefore give trivial fundamental and Brauer groups. In particular $\Pic(X)$ and $H_2(X,\Z)$ are torsion-free, using the exponential and universal coefficient sequences. The same torsion vanishing holds throughout the smooth family. Diaconescu's integral bases modulo torsion \cite[Lemma~2.1]{Diaconescu} now give \eqref{eq:lattices}.

In the chart $z=1$, let $C_\epsilon$ be the section over $D$ given by
\[ x=\epsilon t^6,\qquad y=s^9,\qquad \epsilon\in\{0,\sqrt{-1},-\sqrt{-1}\}. \]
Indeed, $x^3+t^{12}x=0$ for these three values of $x$. Each section is disjoint from $O$, so $r(C_\epsilon)=3$. Any two meet only at $t=0$, where $s=1$ and $(t,x)$ are local coordinates on $Y$ since $2y\ne0$. Their local intersection has length six. By Proposition~\ref{prop:height} the height matrix of $C_0$ and $C_{\sqrt{-1}}$ is $\left(\begin{smallmatrix}6&-3\\-3&6\end{smallmatrix}\right)$, of determinant $27$. Thus $C_0$ and $C_{\sqrt{-1}}$ are independent. The three displayed points are collinear on the generic cubic, so their Mordell--Weil sum is $O$.

Differentiating the Weierstrass equation along $C_\epsilon$ gives
\[ 0\longrightarrow N_{C_\epsilon/X}\longrightarrow \OO_D(1)\oplus\OO_D(6)\oplus\OO_D(9)\xrightarrow{\ d\ }\OO_D(18)\longrightarrow0, \]
where $d=(-s^7t^{10}-\epsilon t^6 A_{11}|_D,-(3\epsilon^2+1)t^{12},2s^9)$. At $A_{11}|_D=0$, the three summands on global sections span the monomials $s^jt^{18-j}$ with $j=7,8$, $j=0,\ldots,6$, and $j=9,\ldots,18$, respectively. Hence $H^0(d)$ is an isomorphism for all three sections for general $A_{11}$. Each normal bundle has degree $-2$ and no global sections, so $N_{C_\epsilon/X}\cong\OO_D(-1)^{\oplus2}$. Inversion gives three further rigid sections with $y=-s^9$. Thus in the single ambient class $\sigma_*[D]+3f$, Theorem~\ref{thm:main} gives at least six isolated points, one for each of the six sections.

More generally, write $[m]C_\epsilon$ for the Mordell--Weil multiple by a nonzero integer $m$. On the smooth locus of $\pi$, the Weierstrass group law is regular \cite[Corollary~12]{Jin}. Multiplication by $m$ is \'etale there since its relative differential is multiplication by $m$ after translation. It maps $C_\epsilon$ isomorphically to $[m]C_\epsilon$ over $D$, identifying their normal bundles. Since $r([m]C_\epsilon)=3m^2$, every such section gives a component isomorphic to $Y^{[n]}$ in class $\sigma_*[D]+(3m^2+n)f$ for each $n\geq0$. For $n=0,1,2$, its GV contributions are
\[ 1,\qquad 36-2w,\qquad 702-78w+3w^2. \]
Let $C=C_0$. At $w=0$, \eqref{eq:main-product} gives $\sum_{n\geq0}n^C_{0,n}q^n=\prod_{m\geq1}(1-q^m)^{-36}$. This is the factor for adding fibers to an isolated section over a line in Greer's genus-zero GW calculation \cite[p.~3, before Conjecture~1.6]{GreerWeierstrass}. Greer uses Noether--Lefschetz theory to count rational curves of base degree one on very general Weierstrass fibrations \cite{GreerNL}.
\end{example}

\begin{corollary}\label{cor:rank-two}
In Example~\ref{ex:plane-line}, put $P=C_0$ and $Q=C_{\sqrt{-1}}$. For very general $A_{11},B_{16}$, every nonzero section in $\Z P+\Z Q$ is infinitesimally rigid in $X$. With $q$ recording total fiber degree, the combined GV contribution of these sections is
\[
\left(\sum_{(a,b)\in\Z^2\setminus\{(0,0)\}}q^{3(a^2-ab+b^2)}\right)\sum_{n\geq0}Q_n(y)q^n,
\]
where $Q_n(y)$ is the contribution from $n$ added fibers in \eqref{eq:hilb-product}.
\end{corollary}
\begin{proof}
Keep $Y$ fixed and specialize $A_{11}|_D=\sqrt{-1}\,s^7t^4$. The coordinates $(\sqrt{-1}\,t^6,s^9)$ then satisfy the Weierstrass equation modulo $u^2$, so $Q$ extends to a section over the doubled line $2D=(u^2=0)$. The normal map for $P$ is independent of $A_{11},B_{16}$, hence $[a]P$ is rigid for $a\ne0$ by Example~\ref{ex:plane-line}. Translation by the extended section $[b]Q$ on the smooth group locus over $2D$ identifies the first infinitesimal neighborhoods of $[a]P$ and $[a]P+[b]Q$ in $X$. Indeed, the square of the ideal of $D$ pulls back into the square of the ideal of either section. Their normal bundles are therefore isomorphic.

For each fixed $(a,b)$ the curve $[a]P+[b]Q\subset Y$ is independent of the coefficients. By upper semicontinuity of $h^0$ of its normal bundle, rigidity holds on a nonempty open set when $a\ne0$. For $a=0$ and $b\ne0$ the same conclusion follows from Example~\ref{ex:plane-line}. Outside the countable union of the complementary closed sets all these sections are rigid. Their half-heights are $3(a^2-ab+b^2)$ by the height matrix above. Theorem~\ref{thm:main}, with $g=0$, gives the stated sum.
\end{proof}

\section{The complete formula over a \texorpdfstring{$(-1)$}{(-1)}-curve}\label{sec:exceptional}

Suppose now that $D=E\cong\PP^1$ with $E^2=-1$ and $Y=\pi^{-1}(E)$ smooth, so that $\chi_Y=1$ and $Y$ is a rational elliptic surface. Let
\[ E_8:=\left\{v\in\Z^8\cup(\Z+\tfrac12)^8\;\middle|\;\sum_{i=1}^8v_i\in2\Z\right\} \]
be the $E_8$ lattice with its standard Euclidean inner product, and let
\begin{equation}\label{eq:theta} \Theta_{E_8}(q):=\sum_{v\in E_8}q^{(v,v)/2} =E_4(q)=1+240\sum_{m\geq1}\Bigl(\sum_{d\mid m}d^3\Bigr)q^m, \end{equation}
be its theta series, where $E_4$ is the normalized Eisenstein series of weight four.

\begin{corollary}\label{cor:E8}
For $\beta_k=\sigma_*[E]+kf$, the full moduli scheme and its canonically oriented sheaf of vanishing cycles are
\begin{equation}\label{eq:E8moduli} M_{\beta_k}\cong\coprod_{\substack{v\in E_8\\(v,v)/2\leq k}} Y^{[k-(v,v)/2]},\qquad \phi_{M_{\beta_k}}\big|_{Y^{[n]}}=\Q_{Y^{[n]}}[2n]. \end{equation}
Here $n=k-(v,v)/2$ on the component indexed by $v$. With $w=y+2+y^{-1}$, the GV invariants satisfy
\begin{equation}\label{eq:E8product} \sum_{k,h\geq0}n_{h,\beta_k}w^hq^k =\Theta_{E_8}(q)\prod_{m\geq1} \frac{1}{(1-q^m)^8(1+yq^m)^2(1+y^{-1}q^m)^2}. \end{equation}
In particular,
\begin{equation}\label{eq:bound} n_{h,\beta_k}=0\quad(h>k),\qquad n_{k,\beta_k}=(-1)^k(k+1), \end{equation}
and the genus-zero series is
\begin{equation}\label{eq:genus-zero} \sum_{k\geq0}n_{0,\beta_k}q^k =\Theta_{E_8}(q)\prod_{m\geq1}(1-q^m)^{-12} =q^{1/2}\frac{E_4(q)}{\eta(q)^{12}}, \qquad\eta(q):=q^{1/24}\prod_{m\geq1}(1-q^m). \end{equation}
\end{corollary}
\begin{proof}
Since $E^2=-1$, any effective divisor on $B$ linearly equivalent to $E$ must contain $E$ with trivial residual divisor, so $|E|$ is a point and $M_k(E)=M_{\beta_k}$. Moreover, for every section $C$, the sequence \eqref{eq:section-normal} becomes
\[ 0\longrightarrow\OO_{\PP^1}(-1)\longrightarrow N_{C/X} \longrightarrow\OO_{\PP^1}(-1)\longrightarrow0, \]
which splits and gives $N_{C/X}\cong\OO_{\PP^1}(-1)^{\oplus2}$, so Theorem~\ref{thm:main} applies to every component.

Because $Y$ is a rational elliptic surface with integral fibers, the orthogonal complement of $O,F$ in its unimodular N\'eron--Severi lattice is negative definite, even (since $K_Y=-F$), and unimodular of rank eight, so its Mordell--Weil lattice is isomorphic to $E_8$ \cite[\S10]{Shioda}. Proposition~\ref{prop:height} therefore gives $r(C)=(v,v)/2$, and combining Theorem~\ref{thm:fixed} with \eqref{eq:main-product} yields \eqref{eq:E8moduli} and \eqref{eq:E8product}.

At degree $q^k$ in \eqref{eq:E8product}, the maximal exponent of $y$ is $k$, which arises solely from the factor $(1+yq)^{-2}$ paired with the constant terms of all other factors and has coefficient $(-1)^k(k+1)$. Since $w^k$ has leading term $y^k$, this gives \eqref{eq:bound}, while setting $w=0$ (that is, $y=-1$) yields \eqref{eq:genus-zero}.
\end{proof}

At $k=1$, the $240$ roots of $E_8$ give the decomposition
\[ M_{\beta_1}\cong Y\sqcup\{240\text{ reduced points}\},\qquad 252-2w=(12-2w)+240, \]
which accounts for the $240$ additional section contributions identified by HKKW for the exceptional curve of the Hirzebruch surface $\mathbb F_1=\PP_{\PP^1}(\OO\oplus\OO(1))$ \cite[\S5.1.4]{HKKW}. The first coefficients of \eqref{eq:E8product} are
\begin{equation}\label{eq:table} \begin{array}{c|l} k&\displaystyle\sum_{h\geq0}n_{h,\beta_k}w^h\\ \hline 0&1\\ 1&252-2w\\ 2&5130-510w+3w^2\\ 3&54760-11780w+772w^2-4w^3\\ 4&419895-142330w+19467w^2-1038w^3+5w^4. \end{array} \end{equation}

The factor multiplying $\Theta_{E_8}$ in \eqref{eq:E8product} is the fixed-section product of Hosono, Saito, and Takahashi \cite[Theorems~4.7 and~4.10]{HST}. Summing over sections supplies the theta series, as in their earlier formula \cite[equation~(5.18)]{HST99}. At genus zero, Bryan and Leung compute the contribution of a fixed section class on the rational elliptic surface \cite[Theorem~6.2]{BryanLeung}. Rose and Yui sum over the $E_8$ lattice in their compact example \cite[Theorems~7.3--7.4]{RoseYui}. Dividing \eqref{eq:E8product} by $w$ also recovers Oberdieck and Shen's PT series in base class $E$, normalized by the pure-fiber PT series \cite[Theorem~2]{OS}.

Klemm, Mayr, and Vafa obtained \eqref{eq:genus-zero} for a singly wound string \cite[equation~(4.1) and Table~1]{KMV}.

At all genera, let $Z_1$ be the elliptic genus of one E-string with all $E_8$ mass parameters zero and rotation parameters $(z,-z)$, as in \cite[equation~(6.30)]{HLV}. With $q=e^{2\pi i\tau}$ and $y=-e^{2\pi iz}$, the Jacobi triple product gives
\[
\sum_{k,h\geq0}n_{h,\beta_k}w^hq^k=q^{1/2}wZ_1,\qquad w=4\sin^2(\pi z).
\]
The compact moduli spaces realizing these coefficients are given in \eqref{eq:E8moduli}.

When $D$ moves, Proposition~\ref{prop:tangent} identifies the extra tangent directions at the fixed-divisor locus with $H^0(C,N_{C/X})$.

\appendix
\section{The direct image for a family of sections}\label{app:family}

Let $U\subset |D|$ parametrize smooth connected curves $D'$ for which $\pi^{-1}(D')$ is smooth. Let $\mathcal D\to U$ be the universal curve and $\mathcal Y:=X\times_B\mathcal D$. Write $\mathcal D^{(n)}_U$ for the relative symmetric product and $\mathcal Y^{[n]}_U$ for the relative Hilbert scheme of $n$ points. The maps $\rho_n$ give the proper morphism
\[ \rho_{U,n}:\mathcal Y^{[n]}_U\longrightarrow\mathcal D^{(n)}_U,\qquad (D',Z)\longmapsto (D',\pi_*[Z]). \]

For each integer $r\geq0$, let $\mathcal S_r$ be the reduced locus in the relative Hilbert scheme of $\mathcal Y/U$ parametrizing pairs $(D',C)$ with $C$ a section satisfying $\langle C,C\rangle=2r$. For $n\geq0$, adding $n$ fibers gives a cycle in class $\beta_{r+n}=\sigma_*[D]+(r+n)f$. Define
\[ \tau_{r,n}:\mathcal S_r\times_U\mathcal D^{(n)}_U\longrightarrow\Chow_{\beta_{r+n}}(X),\qquad (D',C,\textstyle\sum_xm_x[x])\longmapsto[C]+\textstyle\sum_xm_x[f_x]. \]
At $n=0$, $\tau_{r,0}(D',C)=[C]$. Over $U$, write $\mathcal K_k:=Rh_{\beta_k,*}\phi_{M_{\beta_k}}$ for the direct image of vanishing cycles under the support map.

\begin{conjecture}\label{conj:family}
For a choice of orientations that restrict to the canonical ones in Theorem~\ref{thm:main}, there are isomorphisms of constructible complexes on $\mathcal S_r\times_U\mathcal D^{(n)}_U$
\begin{equation}\label{eq:family} \tau_{r,n}^*\mathcal K_{r+n}\cong \operatorname{pr}_1^*(\tau_{r,0}^*\mathcal K_r)\otimes\operatorname{pr}_2^*R\rho_{U,n*}\Q_{\mathcal Y^{[n]}_U}[2n], \end{equation}
where $\operatorname{pr}_1$ and $\operatorname{pr}_2$ are the projections to $\mathcal S_r$ and $\mathcal D^{(n)}_U$, respectively.
\end{conjecture}

To compare with Theorem~\ref{thm:main}, fix $(D,C)\in\mathcal S_r$. On $\{(D,C)\}\times D^{(n)}$, proper base change identifies the right-hand side of \eqref{eq:family} with
\[ (\mathcal K_r)_{[C]}\otimes R\rho_{n*}\Q_{Y^{[n]}}[2n]. \]
If $H^0(C,N_{C/X})=0$, Theorem~\ref{thm:main} at $n=0$ identifies $(\mathcal K_r)_{[C]}$ with $R\Gamma(\Pic^{g+r}(D),\Q)[g]$. For every $n$, the sheaves over these cycles form the open and closed component $M_{C,n}$. On this slice $\tau_{r,n}$ is the finite injective map $\tau_{C,n}$ of Section~\ref{sec:direct}. Pulling back \eqref{eq:fixed-direct} gives \eqref{eq:family}. Taking the perverse Euler polynomial on $D^{(n)}$ gives $w^gQ_n(y)$, the product of the Picard and Hilbert scheme contributions. Summing over $n$ gives \eqref{eq:main-product}. For a moving or obstructed section, \eqref{eq:fixed-direct} computes the intrinsic contribution over $D$, while \eqref{eq:family} retains the ambient stalk $(\mathcal K_r)_{[C]}$.

For $B=\PP^2$ and $D$ of degree $d$, a section $C\ne O$ with $n$ added fibers has total fiber degree $k=r(C)+n$. Its cycle meets $S$ in degree
\[ S\cdot([C]+nf)=C\cdot O+n=k-3d. \]
For the section over a line in Example~\ref{ex:plane-line}, $r(C)=3$. At $n=0$ the curve is disjoint from $S$, as in HKKW's charge-zero contribution \cite[\S5.1.7]{HKKW}. HKKW's Wilson-loop sum ranges over residual plane degrees \cite[\S3.4.1, equation~(3.58)]{HKKW}. Conjecture~\ref{conj:family} is the author's proposed description of the direct image over the smooth-curve locus. Deriving the coefficients of HKKW's sum from that direct image remains open.

\section*{Use of generative AI}
The author used ChatGPT (OpenAI) to help draft and revise the manuscript and to check steps in the proofs. The author takes responsibility for the final text and its mathematical content.

\section*{Acknowledgements}
I thank Sheldon Katz for introducing me to this problem. I also thank Min-xin Huang, Kaiwen Sun, Xin Wang, and Yukinobu Toda for useful discussions. This work was supported by JSPS KAKENHI Grant Number JP25K17226. Kavli IPMU is supported by the World Premier International Research Center Initiative (WPI), MEXT, Japan.

\end{document}